\documentclass{amsart}

\newtheorem{theorem}{Theorem}[section]
\newtheorem{lemma}[theorem]{Lemma}

\theoremstyle{definition}

\theoremstyle{remark}

\numberwithin{equation}{section}

\begin{document}

\title{Necessary and Sufficient Condition for zeros of Derivative of Meromorphic and Entire Functions}

\author{ZhaoKun Ma}
\address{}
\curraddr{YanZhou College, ShanDong Radio and TV University, YanZhou, ShanDong 272100 China}
\email{mzk200408@foxmail.com}
\thanks{}

\author{Lande Ma}
\curraddr{School of Mathematical Sciences, Tongji University, Shanghai, 200092, China}
\email{dzy200408@qq.com}
\thanks{}


\date{March 25st, 2022.}

\dedicatory{}


\begin{abstract}
The main result of this paper shows a totally new necessary and sufficient condition
to determine both real and complex zeros of derivative of all entire and meromorphic functions of one complex variable
in the extended complex plane. By using the theorem, we reprove some results about zeros of derivative of Xi function,
Gamma function and digamma function in a new way.
\end{abstract}

\maketitle

\section*{Introduction}
The complex function of one complex variable begins from Cauchy's times has developed for more than two centuries\cite{Fischer}, and the study of zeros of meromorphic functions together with their derivatives has started about seventy years ago\cite{Hayman1,Edrei2,Ostrovskii}. There are some useful method to study meromorphic and entire functions has been developed by mathematicians, such as the Nevanlinna Theory\cite{Edrei1,Hayman2,Ki}. Many theorems about the zeros of derivatives of meromorphic functions, since the zeros of derivatives of meromorphic functions play an important role in the study of value distributions of meromorphic functions. The recent research about derivatives of meromorphic functions can be found from\cite{Langley1,Langley2,Konstantin,Yamanoi,Xuecheng}.

Let $s=\sigma+it$ be an arbitrary point in the extended complex plane $\mathbb{C}\cup\{\infty\}$\cite{Fischer,Stein},
and let $W(s)=u(\sigma,t)+iv(\sigma,t)$ be a meromorphic function of one complex variable in $\mathbb{C}\cup\{\infty\}$,
by far, the only necessary and sufficient condition to determine if $s$ is a zero of $W^{'}(s)$ is the defination $\frac{dW(s)}{ds}=0$\cite{Lars,Stein}.

Let the point $s_0=\sigma_0+it_0\in\mathbb{C}\cup\{\infty\}$ not zero or pole of $W(s)$, then $s_0$ is the zero of $W^{'}(s)$ if and only if two partial derivatives equations $\frac{\partial{\varphi}}{\partial{\sigma}}(\sigma,t)\mid{_{s=s_0}}=0$
and $\frac{\partial{\varphi}}{\partial{t}}(\sigma,t)\mid{_{s=s_0}}=0$ hold simultaneously.
In which, $\varphi(\sigma,t)=arg(\frac{v(\sigma,t)}{u(\sigma,t)})$ is the argument function of $W(s)$.
The main theorem in this paper needn't other requirements and has no other limitation, like the degree of the functions. Except the definition, the main theorem in this paper is the second result of the necessary and sufficient condition.

Because $\varphi(\sigma,t)$ is a real value function which contains two real variables $\sigma$ and $t$.
Its two partial derivatives $\frac{\partial{\varphi}}{\partial{\sigma}}(\sigma,t)$ and $\frac{\partial{\varphi}}{\partial{t}}(\sigma,t)$ are also two real values functions of two real variables. In the theorem which we have shown, the zeros of derivative of meromorphic functions are the zeros of two real partial derivatives.

By utilizing two partial derivatives of arguments, we prove the main theorem,
it is a necessary and sufficient condition to determine both real and complex zeros of derivative of
all the meromorphic functions of one complex variable in $\mathbb{C}\cup\{\infty\}$. The main theorem in this paper can let the problem of the complex variable function of one variable utilize the real variable function tool to solve. The argument of products of factors is equal to the sum of arguments of factors,
it changes the multiplication relation into the addition relation. For the functions that with factor multiplication, their zeros of derivative can be easily obtained by using the new method. Here, in order to prove the valuable of the main theorem, utilizing the main theorem, we solve two problems.

\section{The results}
Let $W(s)=u(\sigma,t)+iv(\sigma,t)$ be a meromorphic function with one complex variable in the extended complex plane $\mathbb{C}\cup\{\infty\}$, we show the expressions of its argument\cite{Fischer,Stein}.

(1). When $u(\sigma,t)>0$, $v(\sigma,t)\ge0$. $\varphi(\sigma,t)=arg(\frac{v(\sigma,t)}{u(\sigma,t)})=\arctan(\frac{v(\sigma,t)}{u(\sigma,t)})$.

(2). When $u(\sigma,t)<0$, $v(\sigma,t)\ge0$. $\varphi(\sigma,t)=arg(\frac{v(\sigma,t)}{u(\sigma,t)})=\pi-\arctan(\frac{v(\sigma,t)}{-u(\sigma,t)})$.

(3). When $u(\sigma,t)<0$, $v(\sigma,t)<0$. $\varphi(\sigma,t)=arg(\frac{v(\sigma,t)}{u(\sigma,t)})=\pi+\arctan(\frac{v(\sigma,t)}{u(\sigma,t)})$.

(4). When $u(\sigma,t)>0$, $v(\sigma,t)<0$. $\varphi(\sigma,t)=arg(\frac{v(\sigma,t)}{u(\sigma,t)})=-\arctan(\frac{-v(\sigma,t)}{u(\sigma,t)})$.

At zeros of $W(s)$, we have $W(s)=0$, the argument $\varphi(\sigma,t)$ can not be surely determined. At poles of $W(s)$, we have $W(s)=\infty$, the argument $\varphi(\sigma,t)$ can not be surely determined either.

For the points which let $u(\sigma,t)=0$,
they can not be in the definition domain of the fraction $\frac{v(\sigma,t)}{u(\sigma,t)}$.
The function $\arctan$ is continuous and differential on the two-dimensional plane formed by the variables $\sigma$ and $t$\cite{Stein}.

The argument function $\varphi(\sigma,t)$ of $W(s)$ is a continuous and differential function of two variables $\sigma$ and $t$,
except zeros and poles and the points which let $u(\sigma,t)$ of $W(s)$ equals to $0$.

\begin{lemma}
Let $W^{'}(s)$ be the derivative of $W(s)$,
and let $\frac{\partial{W}}{\partial{\sigma}}(\sigma+it)$ be the partial derivative of $W(s)$ concerning the variable $\sigma$,
then the zeros of $W^{'}(s)$ and zeros of $\frac{\partial{W}}{\partial{\sigma}}(\sigma+it)$
are all the same points in $\mathbb{C}\cup\{\infty\}$.
\end{lemma}
\begin{proof}
Let $W(s)=\frac{W_1(s)}{W_2(s)}$, its derivative concerning the complex variable $s$ is:
$W^{'}(s)=\frac{W^{'}_1(s)W_2(s)-W^{'}_2(s)W_1(s)}{W^2_2(s)}$.

Substitute $\sigma+it$ for $s$ in $W(s)$, we have: $W(\sigma+it)$.
The $W(\sigma+it)$ can be considered as the function $W(s)$ composed with $\sigma+it$.

We have:
$\frac{\partial{W}}{\partial{\sigma}}(\sigma+it)=
\frac{\frac{dW_1(\sigma+it)}{d(\sigma+it)}\frac{\partial{(\sigma+it)}}{\partial{\sigma}}W_2(\sigma+it)-
\frac{dW_2(\sigma+it)}{d(\sigma+it)}\frac{\partial{(\sigma+it)}}{\partial{\sigma}}W_1(\sigma+it)}{W^2_2(\sigma+it)}$

$=\frac{\frac{dW_1(\sigma+it)}{d(\sigma+it)}W_2(\sigma+it)-
\frac{dW_2(\sigma+it)}{d(\sigma+it)}W_1(\sigma+it)}{W^2_2(\sigma+it)}$.

In the front derivative, let $\sigma+it=s$, we have:

$\frac{\partial{W}}{\partial{\sigma}}(s)=
\frac{\frac{dW_1(s)}{ds}W_2(s)-\frac{dW_2(s)}{ds}W_1(s)}{W^2_2(s)}=\frac{W^{'}_1(s)W_2(s)-W^{'}_2(s)W_1(s)}{W^2_2(s)}$.
The two derivatives $W^{'}(s)$ and $\frac{\partial{W}}{\partial{\sigma}}(s)$ are completely the same.
So the zeros of $W^{'}(s)$ and $\frac{\partial{W}}{\partial{\sigma}}(s)$ are completely the same.
\end{proof}
Applying the same method as Lemma 1.1, we can prove Lemma 1.2.

\begin{lemma}
Let $W^{'}(s)$ be the derivative of $W(s)$,
and let $\frac{\partial{W}}{\partial{t}}(\sigma+it)$ be the partial derivative of $W(s)$ concerning the variable $t$,
then the zeros of $W^{'}(s)$ and zeros of $\frac{\partial{W}}{\partial{t}}(\sigma+it)$ are all the same points in
$\mathbb{C}\cup\{\infty\}$.
\end{lemma}

Let $\frac{\partial{\varphi}}{\partial{\sigma}}(\sigma,t)$
be the partial derivative of argument of $W(s)$ concerning the variable $\sigma$,
$\frac{\partial{\varphi}}{\partial{t}}(\sigma,t)$
be the partial derivative of argument of $W(s)$ concerning the variable $t$.

Let $\Phi=\{s|s=\sigma+it\in\mathbb{C}\cup\{\infty\}$, $s$ is not zero or pole of $W(s)$,
and there is $u(\sigma,t)\neq0\}$.

\begin{lemma}
Let $s=(\sigma,t)\in\Phi$, then we have two partial derivatives:
$\frac{\partial{\varphi}}{\partial{\sigma}}(\sigma,t)=\frac{v^{'}_{\sigma}(\sigma,t)u(\sigma,t)-u^{'}_{\sigma}(\sigma,t)v(\sigma,t)}{u^2(\sigma,t)+v^2(\sigma,t)}$,
$\frac{\partial{\varphi}}{\partial{t}}(\sigma,t)=\frac{v^{'}_{t}(\sigma,t)u(\sigma,t)-u^{'}_{t}(\sigma,t)v(\sigma,t)}{u^2(\sigma,t)+v^2(\sigma,t)}$.
\end{lemma}
\begin{proof}
According to the positive or negative of $u(\sigma,t)$ and $v(\sigma,t)$, we divide our proof into four different cases.

(1). When $u(\sigma,t)>0$, $v(\sigma,t)\ge0$.

$\frac{\partial{\varphi}}{\partial{\sigma}}(\sigma,t)=
\frac{\partial{\arctan(\frac{v(\sigma,t)}{u(\sigma,t)})}}{\partial{\sigma}}
=\frac{1}{1+(\frac{v(\sigma,t)}{u(\sigma,t)})^2}\frac{v^{'}_{\sigma}(\sigma,t)u(\sigma,t)-u^{'}_{\sigma}(\sigma,t)v(\sigma,t)}{u^2(\sigma,t)}$

$=\frac{v^{'}_{\sigma}(\sigma,t)u(\sigma,t)-u^{'}_{\sigma}(\sigma,t)v(\sigma,t)}{u^2(\sigma,t)+v^2(\sigma,t)}$,

$\frac{\partial{\varphi}}{\partial{t}}(\sigma,t)=
\frac{\partial{\arctan(\frac{v(\sigma,t)}{u(\sigma,t)})}}{\partial{t}}
=\frac{1}{1+(\frac{v(\sigma,t)}{u(\sigma,t)})^2}\frac{v^{'}_{t}(\sigma,t)u(\sigma,t)-u^{'}_{t}(\sigma,t)v(\sigma,t)}{u^2(\sigma,t)}$

$=\frac{v^{'}_{t}(\sigma,t)u(\sigma,t)-u^{'}_{t}(\sigma,t)v(\sigma,t)}{u^2(\sigma,t)+v^2(\sigma,t)}$.

(2). When $u(\sigma,t)<0$, $v(\sigma,t)\ge0$.

$\frac{\partial{\varphi}}{\partial{\sigma}}(\sigma,t)=
\frac{\partial{(\pi-\arctan(\frac{v(\sigma,t)}{-u(\sigma,t)}))}}{\partial{\sigma}}
=-\frac{1}{1+(\frac{v(\sigma,t)}{u(\sigma,t)})^2}\frac{(-1)(v^{'}_{\sigma}(\sigma,t)u(\sigma,t)-u^{'}_{\sigma}(\sigma,t)v(\sigma,t))}{u^2(\sigma,t)}$

$=\frac{v^{'}_{\sigma}(\sigma,t)u(\sigma,t)-u^{'}_{\sigma}(\sigma,t)v(\sigma,t)}{u^2(\sigma,t)+v^2(\sigma,t)}$,

$\frac{\partial{\varphi}}{\partial{t}}(\sigma,t)=
\frac{\partial{(\pi-\arctan(\frac{v(\sigma,t)}{-u(\sigma,t)}))}}{\partial{t}}
=-\frac{1}{1+(\frac{v(\sigma,t)}{u(\sigma,t)})^2}\frac{(-1)(v^{'}_{t}(\sigma,t)u(\sigma,t)-u^{'}_{t}(\sigma,t)v(\sigma,t))}{u^2(\sigma,t)}$

$=\frac{v^{'}_{t}(\sigma,t)u(\sigma,t)-u^{'}_{t}(\sigma,t)v(\sigma,t)}{u^2(\sigma,t)+v^2(\sigma,t)}$.

(3). When $u(\sigma,t)<0$, $v(\sigma,t)<0$.

$\frac{\partial{\varphi}}{\partial{\sigma}}(\sigma,t)=
\frac{\partial{(\pi+\arctan(\frac{v(\sigma,t)}{u(\sigma,t)}))}}{\partial{\sigma}}
=\frac{1}{1+(\frac{v(\sigma,t)}{u(\sigma,t)})^2}\frac{v^{'}_{\sigma}(\sigma,t)u(\sigma,t)-u^{'}_{\sigma}(\sigma,t)v(\sigma,t)}{u^2(\sigma,t)}$

$=\frac{v^{'}_{\sigma}(\sigma,t)u(\sigma,t)-u^{'}_{\sigma}(\sigma,t)v(\sigma,t)}{u^2(\sigma,t)+v^2(\sigma,t)}$,

$\frac{\partial{\varphi}}{\partial{t}}(\sigma,t)=
\frac{\partial{(\pi+\arctan(\frac{v(\sigma,t)}{u(\sigma,t)}))}}{\partial{t}}
=\frac{1}{1+(\frac{v(\sigma,t)}{u(\sigma,t)})^2}\frac{v^{'}_{t}(\sigma,t)u(\sigma,t)-u^{'}_{t}(\sigma,t)v(\sigma,t)}{u^2(\sigma,t)}$

$=\frac{v^{'}_{t}(\sigma,t)u(\sigma,t)-u^{'}_{t}(\sigma,t)v(\sigma,t)}{u^2(\sigma,t)+v^2(\sigma,t)}$.

(4). When $u(\sigma,t)>0$, $v(\sigma,t)<0$.

$\frac{\partial{\varphi}}{\partial{\sigma}}(\sigma,t)=
\frac{\partial{(-\arctan(\frac{-v(\sigma,t)}{u(\sigma,t)}))}}{\partial{\sigma}}
=-\frac{1}{1+(\frac{v(\sigma,t)}{u(\sigma,t)})^2}\frac{(-1)(v^{'}_{\sigma}(\sigma,t)u(\sigma,t)-u^{'}_{\sigma}(\sigma,t)v(\sigma,t))}{u^2(\sigma,t)}$

$=\frac{v^{'}_{\sigma}(\sigma,t)u(\sigma,t)-u^{'}_{\sigma}(\sigma,t)v(\sigma,t)}{u^2(\sigma,t)+v^2(\sigma,t)}$,

$\frac{\partial{\varphi}}{\partial{t}}(\sigma,t)=
\frac{\partial{(-\arctan(\frac{-v(\sigma,t)}{u(\sigma,t)}))}}{\partial{t}}
=-\frac{1}{1+(\frac{v(\sigma,t)}{u(\sigma,t)})^2}\frac{(-1)(v^{'}_{t}(\sigma,t)u(\sigma,t)-u^{'}_{t}(\sigma,t)v(\sigma,t))}{u^2(\sigma,t)}$

$=\frac{v^{'}_{t}(\sigma,t)u(\sigma,t)-u^{'}_{t}(\sigma,t)v(\sigma,t)}{u^2(\sigma,t)+v^2(\sigma,t)}$.

Sum up the above four cases, we have two partial derivatives:

$\frac{\partial{\varphi}}{\partial{\sigma}}(\sigma,t)=
\frac{v^{'}_{\sigma}(\sigma,t)u(\sigma,t)-u^{'}_{\sigma}(\sigma,t)v(\sigma,t)}{u^2(\sigma,t)+v^2(\sigma,t)}$,
$\frac{\partial{\varphi}}{\partial{t}}(\sigma,t)=
\frac{v^{'}_{t}(\sigma,t)u(\sigma,t)-u^{'}_{t}(\sigma,t)v(\sigma,t)}{u^2(\sigma,t)+v^2(\sigma,t)}$
\end{proof}

\begin{lemma}
Let $s_0=\sigma_0+it_0\in\Phi$,
then the necessary and sufficient condition which two equations
$\frac{\partial{\varphi}}{\partial{\sigma}}(\sigma,t)\mid{_{s=s_0}}=0$ and $\frac{\partial{\varphi}}{\partial{t}}(\sigma,t)\mid{_{s=s_0}}=0$
hold simultaneously is that two groups of the equations
$\frac{\partial{u}}{\partial{\sigma}}(\sigma,t)\mid{_{s=s_0}}=\frac{\partial{v}}{\partial{\sigma}}(\sigma,t)\mid{_{s=s_0}}=0$ and
$\frac{\partial{u}}{\partial{t}}(\sigma,t)\mid{_{s=s_0}}=\frac{\partial{v}}{\partial{t}}(\sigma,t)\mid{_{s=s_0}}=0$ hold simultaneously.
\end{lemma}
\begin{proof}
Sufficiency. According to the Cauchy-Riemann Equations\cite{Fischer},
since $W(s)$ is a meromorphic function, for any point $s=\sigma+it\in\Phi$, we have $v^{'}_{\sigma}(\sigma,t)=-u^{'}_{t}(\sigma,t)$ and $v^{'}_{t}(\sigma,t)=u^{'}_{\sigma}(\sigma,t)$.

At any point $s_0\in\Phi$, when two equations $\frac{\partial{u}}{\partial{\sigma}}(\sigma,t)\mid{_{s=s_0}}=
\frac{\partial{v}}{\partial{\sigma}}(\sigma,t)\mid{_{s=s_0}}=0$ are true.
$\frac{\partial{u}}{\partial{t}}(\sigma,t)\mid{_{s=s_0}}=
-\frac{\partial{v}}{\partial{\sigma}}(\sigma,t)\mid{_{s=s_0}}$,
$\frac{\partial{u}}{\partial{t}}(\sigma,t)\mid{_{s=s_0}}=0$.
And $\frac{\partial{v}}{\partial{t}}(\sigma,t)\mid{_{s=s_0}}=\frac{\partial{u}}{\partial{\sigma}}(\sigma,t)\mid{_{s=s_0}}$, $\frac{\partial{v}}{\partial{t}}(\sigma,t)\mid{_{s=s_0}}=0$. Namely, two equations $\frac{\partial{u}}{\partial{t}}(\sigma,t)\mid{_{s=s_0}}=\frac{\partial{v}}{\partial{t}}(\sigma,t)\mid{_{s=s_0}}=0$ hold.

Conversely, for any point $s_0\in\Phi$, when two equations $\frac{\partial{u}}{\partial{t}}(\sigma,t)\mid{_{s=s_0}}=\frac{\partial{v}}{\partial{t}}(\sigma,t)\mid{_{s=s_0}}=0$ hold, $\frac{\partial{u}}{\partial{\sigma}}(\sigma,t)\mid{_{s=s_0}}=\frac{\partial{v}}{\partial{t}}(\sigma,t)\mid{_{s=s_0}}$, $\frac{\partial{u}}{\partial{\sigma}}(\sigma,t)\mid{_{s=s_0}}=0$. And $\frac{\partial{v}}{\partial{\sigma}}(\sigma,t)\mid{_{s=s_0}}=-\frac{\partial{u}}{\partial{t}}(\sigma,t)\mid{_{s=s_0}}$, $\frac{\partial{v}}{\partial{\sigma}}(\sigma,t)\mid{_{s=s_0}}=0$. Namely, two equations $\frac{\partial{u}}{\partial{\sigma}}(\sigma,t)\mid{_{s=s_0}}=\frac{\partial{v}}{\partial{\sigma}}(\sigma,t)\mid{_{s=s_0}}=0$ hold.

The above four equations separately multiply two function values $(u(\sigma_0,t_0)$ and $v(\sigma_0,t_0))$ which only have the finite values, the four equations $v^{'}_{\sigma}(\sigma,t)\mid{_{s=s_0}}u(\sigma_0,t_0)=0$, $u^{'}_{\sigma}(\sigma,t)\mid{_{s=s_0}}v(\sigma_0,t_0)=0$. $v^{'}_{t}(\sigma,t)\mid{_{s=s_0}}u(\sigma_0,t_0)=0$, $u^{'}_{t}(\sigma,t)\mid{_{s=s_0}}v(\sigma_0,t_0)=0$ hold simultaneously. $u^2(\sigma_0,t_0)+v^2(\sigma_0,t_0)\neq 0$. According to two partial derivatives in Lemma 1.3, we have: $\frac{\partial{\varphi}}{\partial{\sigma}}(\sigma,t)\mid{_{s=s_0}}=0$ and $\frac{\partial{\varphi}}{\partial{t}}(\sigma,t)\mid{_{s=s_0}}=0$.

Necessity. Two equations $\frac{\partial{\varphi}}{\partial{\sigma}}(\sigma,t)\mid{_{s=s_0}}=0$ and
$\frac{\partial{\varphi}}{\partial{t}}(\sigma,t)\mid{_{s=s_0}}=0$ hold simultaneously,
namely, $\frac{\partial{\varphi}}{\partial{\sigma}}(\sigma,t)\mid{_{s=s_0}}=\frac{v^{'}_{\sigma}(\sigma,t)\mid{_{s=s_0}}u(\sigma_0,t_0)-u^{'}_{\sigma}(\sigma,t)\mid{_{s=s_0}}v(\sigma_0,t_0)}{u^2(\sigma_0,t_0)+v^2(\sigma_0,t_0)}=0$,

$\frac{\partial{\varphi}}{\partial{t}}(\sigma,t)\mid{_{s=s_0}}=
\frac{v^{'}_{t}(\sigma,t)\mid{_{s=s_0}}u(\sigma_0,t_0)-u^{'}_{t}(\sigma,t)\mid{_{s=s_0}}v(\sigma_0,t_0)}{u^2(\sigma_0,t_0)+v^2(\sigma_0,t_0)}=0$.

In this lemma, we have a constraint: the point $s_0$ is in $\Phi$.
For any point in $\Phi$, the function $u^2(\sigma,t)+v^2(\sigma,t)$ obtains the non-zero finite value.
So, if the two equations are true, then, only their numerator polynomials are equal to 0.
So, we have:  $v^{'}_{\sigma}(\sigma_0,t_0)u(\sigma_0,t_0)-u^{'}_{\sigma}(\sigma_0,t_0)v(\sigma_0,t_0)=0$,
and $v^{'}_{t}(\sigma_0,t_0)u(\sigma_0,t_0)-u^{'}_{t}(\sigma_0,t_0)v(\sigma_0,t_0)=0$.
Namely, $v^{'}_{\sigma}(\sigma_0,t_0)u(\sigma_0,t_0)=u^{'}_{\sigma}(\sigma_0,t_0)v(\sigma_0,t_0)$,
and $v^{'}_{t}(\sigma_0,t_0)$ $u(\sigma_0,t_0)$$=u^{'}_{t}(\sigma_0,t_0)$$v(\sigma_0,t_0)$.

According to the Cauchy-Riemann Equations, for any point in $\Phi$,
we have: $v^{'}_{\sigma}(\sigma_0,t_0)=-u^{'}_{t}(\sigma_0,t_0)$,
and $v^{'}_{t}(\sigma_0,t_0)=u^{'}_{\sigma}(\sigma_0,t_0)$.
Together with the two equations we obtained in last paragraph, we have:
$v^{'}_{\sigma}(\sigma_0,t_0)u(\sigma_0,t_0)=v^{'}_{t}(\sigma_0,t_0)v(\sigma_0,t_0)$ and
$v^{'}_{t}(\sigma_0,t_0)u(\sigma_0,t_0)=-v^{'}_{\sigma}(\sigma_0,t_0)v(\sigma_0,t_0)$.
The former equation multiplies $v(\sigma_0,t_0)$ on both sides,
and the latter equation multiplies $u(\sigma_0,t_0)$ on both sides,
then we have: $v^{'}_{\sigma}(\sigma_0,t_0)u(\sigma_0,t_0)v(\sigma_0,t_0)=v^{'}_{t}(\sigma_0,t_0)v^2(\sigma_0,t_0)$
and $v^{'}_{t}(\sigma_0,t_0)u^2(\sigma_0,t_0)=-v^{'}_{\sigma}(\sigma_0,t_0)v(\sigma_0,t_0)u(\sigma_0,t_0)$.
Two sides of two equations separately plus, we have: $(u^2(\sigma_0,t_0)+v^2(\sigma_0,t_0))v^{'}_{t}(\sigma_0,t_0)=0$.

For any point in $\Phi$, the function $u^2(\sigma,t)+v^2(\sigma,t)$ obtains the non-zero finite value.
$(u^2(\sigma_0,t_0)+v^2(\sigma_0,t_0))\neq 0$, so, $v^{'}_{t}(\sigma_0,t_0)=0$,
then we have: $u^{'}_{\sigma}(\sigma,t)=v^{'}_{t}(\sigma_0,t_0)=0$.

Because $v^{'}_{\sigma}(\sigma_0,t_0)u(\sigma_0,t_0)=u^{'}_{\sigma}(\sigma_0,t_0)v(\sigma_0,t_0)$,
and $v^{'}_{t}(\sigma_0,t_0)u(\sigma_0,t_0)=$ $u^{'}_{t}(\sigma_0,t_0)$ $v(\sigma_0,t_0)$, we have:
$v^{'}_{\sigma}(\sigma_0,t_0)u(\sigma_0,t_0)=0$, $u^{'}_{t}(\sigma_0,t_0)v(\sigma_0,t_0)=0$.
Except zeros of $W(s)$, the two functions $u(\sigma,t)$, $v(\sigma,t)$ cannot simultaneously equal to 0 on the point $s_0=\sigma_0+it_0$. For any point in $\Phi$, two functions $u(\sigma,t)$ and $v(\sigma,t)$ both obtain finite values. So, we have $v^{'}_{\sigma}(\sigma_0,t_0)=0$ and $u^{'}_{t}(\sigma_0,t_0)=0$. And because $v^{'}_{\sigma}(\sigma,t)=-u^{'}_{t}(\sigma,t)$, we have $v^{'}_{\sigma}(\sigma_0,t_0)=u^{'}_{t}(\sigma_0,t_0)=0$.
\end{proof}

\begin{lemma}
Let $s_0\in\Phi$, then $s_0$ is the zero of $W^{'}(s)$
if and only if there are $\frac{\partial{\varphi}}{\partial{\sigma}}(\sigma,t)\mid{_{s=s_0}}=0$ and
$\frac{\partial{\varphi}}{\partial{t}}(\sigma,t)\mid{_{s=s_0}}=0$.
\end{lemma}
\begin{proof}
Sufficiency. If two equations $\frac{\partial{\varphi}}{\partial{\sigma}}(\sigma,t)\mid{_{s=s_0}}=0$ and  $\frac{\partial{\varphi}}{\partial{t}}(\sigma,t)\mid{_{s=s_0}}=0$ hold simultaneously. According to Lemma 1.4, $\frac{\partial{u}}{\partial{\sigma}}(\sigma,t)\mid{_{s=s_0}}=\frac{\partial{v}}{\partial{\sigma}}(\sigma,t)\mid{_{s=s_0}}=0$ and  $\frac{\partial{u}}{\partial{t}}(\sigma,t)\mid{_{s=s_0}}=\frac{\partial{v}}{\partial{t}}(\sigma,t)\mid{_{s=s_0}}=0$ hold simultaneously.

If $\frac{\partial{u}}{\partial{\sigma}}(\sigma,t)\mid{_{s=s_0}}=\frac{\partial{v}}{\partial{\sigma}}(\sigma,t)\mid{_{s=s_0}}=0$, then $\frac{\partial{W}}{\partial{\sigma}}(\sigma,t)=\frac{\partial{u}}{\partial{\sigma}}(\sigma,t)+i\frac{\partial{v}}{\partial{\sigma}}(\sigma,t)=0$,
so, the point $s_0=\sigma_0+it_0\in\Phi$ is a zero of $\frac{\partial{W}}{\partial{\sigma}}(\sigma,t)$. According to Lemma 1.1, the point $s_0=\sigma_0+it_0$ is also a zero of $W^{'}(s)_s$.

If $\frac{\partial{u}}{\partial{t}}(\sigma,t)\mid{_{s=s_0}}=\frac{\partial{v}}{\partial{t}}(\sigma,t)\mid{_{s=s_0}}=0$,
then $\frac{\partial{W}}{\partial{t}}(\sigma,t)=\frac{\partial{u}}{\partial{t}}(\sigma,t)+i\frac{\partial{v}}{\partial{t}}(\sigma,t)=0$, so, the point $s_0=\sigma_0+it_0\in\Phi$ is a zero of $\frac{\partial{W}}{\partial{t}}(\sigma,t)$. According to Lemma 1.2, the point $s_0=\sigma_0+it_0$ is also a zero of $W^{'}(s)$.

Necessity. Because $W(s)=u(\sigma,t)+iv(\sigma,t)$, we do the partial derivative concerning the variable $\sigma$, obtain:
$\frac{\partial{W}}{\partial{\sigma}}(\sigma,t)=\frac{\partial{u}}{\partial{\sigma}}(\sigma,t)+i\frac{\partial{v}}{\partial{\sigma}}(\sigma,t)$.
At the point $s_0$, $\frac{\partial{W}}{\partial{\sigma}}(\sigma,t)\mid{_{s=s_0}}=0$. So, the real and imaginary part of $\frac{\partial{W}}{\partial{\sigma}}(\sigma,t)\mid{_{s=s_0}}$ must be separately 0. Namely, we have  $\frac{\partial{u}}{\partial{\sigma}}(\sigma,t)=\frac{\partial{v}}{\partial{\sigma}}(\sigma,t)=0$. According to Lemma 1.4, for any point in $\Phi$, we have $\frac{\partial{\varphi}}{\partial{\sigma}}(\sigma,t)\mid{_{s=s_0}}=0$ and $\frac{\partial{\varphi}}{\partial{t}}(\sigma,t)\mid{_{s=s_0}}=0$.

According to Lemma 1.1, the zeros of $W^{'}(s)$ and the zeros of $\frac{\partial{W}}{\partial{\sigma}}(\sigma,t)$ are completely the same points in $\mathbb{C}\cup\{\infty\}$.

So, we have: the zero $s_0\in\Phi$ of derivative of $W(s)$ can let two equations $\frac{\partial{\varphi}}{\partial{\sigma}}(\sigma,t)\mid{_{s=s_0}}=0$ and $\frac{\partial{\varphi}}{\partial{t}}(\sigma,t)\mid{_{s=s_0}}=0$ hold simultaneously.

We do the partial derivative concerning the variable $t$ for $W(s)=u(\sigma,t)+iv(\sigma,t)$, obtain:
$\frac{\partial{W}}{\partial{t}}(\sigma,t)=\frac{\partial{u}}{\partial{t}}(\sigma,t)+i\frac{\partial{v}}{\partial{t}}(\sigma,t)$. At the point  $s_0$, $\frac{\partial{W}}{\partial{t}}(\sigma,t)\mid{_{s=s_0}}=0$. So, the real and imaginary part of $\frac{\partial{W}}{\partial{t}}(\sigma,t)\mid{_{s=s_0}}$ must be separately 0. Namely, we have  $\frac{\partial{u}}{\partial{t}}(\sigma,t)=\frac{\partial{v}}{\partial{t}}(\sigma,t)=0$. According to Lemma 1.4, for any point in $\Phi$, we have $\frac{\partial{\varphi}}{\partial{\sigma}}(\sigma,t)\mid{_{s=s_0}}=0$ and  $\frac{\partial{\varphi}}{\partial{t}}(\sigma,t)\mid{_{s=s_0}}=0$.

According to Lemma 1.2, the zeros of $W^{'}(s)$, the zeros
of $\frac{\partial{W}}{\partial{t}}(\sigma,t)$ are totally the same points in $\mathbb{C}\cup\{\infty\}$.

So, we obtain: the zero $s_0\in\Phi$ of $W^{'}(s)$ let two equations $\frac{\partial{\varphi}}{\partial{\sigma}}(\sigma,t)\mid{_{s=s_0}}=0$
and $\frac{\partial{\varphi}}{\partial{t}}(\sigma,t)\mid{_{s=s_0}}=0$ hold simultaneously.
\end{proof}

\begin{lemma}
Let $s_0=\sigma_0+it_0\in\Phi$, if $s_0$ is not zero of $W^{'}(s)$,
then $\frac{\partial{\varphi}}{\partial{\sigma}}(\sigma,t)$
and $\frac{\partial{\varphi}}{\partial{t}}(\sigma,t)$ cannot simultaneously equal to 0.
\end{lemma}

The lemma can also be written as: let $s_0=\sigma_0+it_0\in\Phi$,
if it is not zero of $W^{'}(s)$, then if $\frac{\partial{\varphi}}{\partial{\sigma}}(\sigma,t)\mid{_{s=s_0}}=0$,
we have: $\frac{\partial{\varphi}}{\partial{t}}(\sigma,t)\mid{_{s=s_0}}\neq 0$; or,
if $\frac{\partial{\varphi}}{\partial{t}}(\sigma,t)\mid{_{s=s_0}}=0$, we have:
$\frac{\partial{\varphi}}{\partial{\sigma}}(\sigma,t)\mid{_{s=s_0}}\neq 0$.

\begin{proof}
According to Lemma 1.5, if two equations $\frac{\partial{\varphi}}{\partial{\sigma}}(\sigma,t)\mid{_{s=s_0}}=0$
and $\frac{\partial{\varphi}}{\partial{t}}(\sigma,t)\mid{_{s=s_0}}=0$ hold simultaneously,
then  $s_0$ is surely a zero of $W^{'}(s)$. So, for any point in $\Phi$,
except the zeros of $W^{'}(s)$, the partial derivative of argument of $W(s)$ concerning the variable $\sigma$ and variable $t$ cannot simultaneously equal to 0.
\end{proof}

In the results and their proofs in front, the situation of $u(\sigma,t)=0$ is excluded, so, the following, we need to prove the same results are true under the situation of $u(\sigma,t)=0$.
Let $W_3(s)=iW(s)$, so, $W_3(s)=iu(\sigma,t)-v(\sigma,t)$. It is obvious that $W_3(s)$ is a meromorphic function.
We can apply the results of Lemma 1.1 and Lemma 1.2 to $W_3(s)$,
so, the same results concerning $W_3(s)$ don't need to be given.

And the argument of function $W_3(s)$ is $\varphi_3(\sigma,t)=\frac{\pi}{2}+\varphi(\sigma,t)$.
$\frac{\partial{\varphi_3}}{\partial{\sigma}}(\sigma,t)=\frac{\partial{\varphi}}{\partial{\sigma}}(\sigma,t)$,
$\frac{\partial{\varphi_3}}{\partial{t}}(\sigma,t)=\frac{\partial{\varphi}}{\partial{t}}(\sigma,t)$.
So, we obtain Lemma 1.7.

Let $\Omega=\{s|s=\sigma+it\in\mathbb{C}\cup\{\infty\}$, $s$ is not zero or pole of $W_3(s)$,
and there is $v(\sigma,t)\neq0\}$.

\begin{lemma}
Let $s=\sigma+it\in\Omega$, then we have two partial derivatives:
$\frac{\partial{\varphi_3}}{\partial{\sigma}}(\sigma,t)=
\frac{v^{'}_{\sigma}(\sigma,t)u(\sigma,t)-u^{'}_{\sigma}(\sigma,t)v(\sigma,t)}{u^2(\sigma,t)+v^2(\sigma,t)}$,
$\frac{\partial{\varphi_3}}{\partial{t}}(\sigma,t)=
\frac{v^{'}_{t}(\sigma,t)u(\sigma,t)-u^{'}_{t}(\sigma,t)v(\sigma,t)}{u^2(\sigma,t)+v^2(\sigma,t)}$
\end{lemma}

According to the Cauchy-Riemann Equations, and $W_3(s)$ is a meromorphic function, for any point in $\Omega$,
we have: $v^{'}_{\sigma}(\sigma,t)=-u^{'}_{t}(\sigma,t)$, $v^{'}_{t}(\sigma,t)=u^{'}_{\sigma}(\sigma,t)$.
So, the conclusions of Lemma 1.4, Lemma 1.5 and Lemma 1.6 can apply to $\varphi_3(\sigma,t)$ and its derivatives.
Namely, the proofs of Lemma 1.8, Lemma 1.9 and Lemma 1.10 are the same as the proofs of the three corresponding lemmas in front.

\begin{lemma}
Let $s_0=\sigma_0+it_0\in\Omega$, then the necessary and sufficient condition which two equations
$\frac{\partial{\varphi_3}}{\partial{\sigma}}(\sigma,t)\mid{_{s=s_0}}=0$
and $\frac{\partial{\varphi_3}}{\partial{t}}(\sigma,t)\mid{_{s=s_0}}=0$
hold simultaneously is that two groups of the equations
$\frac{\partial{u}}{\partial{\sigma}}(\sigma,t)\mid{_{s=s_0}}=
\frac{\partial{v}}{\partial{\sigma}}(\sigma,t)\mid{_{s=s_0}}=0$ and the equations
$\frac{\partial{u}}{\partial{t}}(\sigma,t)\mid{_{s=s_0}}=
\frac{\partial{v}}{\partial{t}}(\sigma,t)\mid{_{s=s_0}}=0$ hold simultaneously.
\end{lemma}

\begin{lemma}
Let $s_0\in\Omega$, then the point $s_0$ is the zero of $W_3^{'}(s)$
if and only if there are two equations $\frac{\partial{\varphi_3}}{\partial{\sigma}}(\sigma,t)\mid{_{s=s_0}}=0$
and $\frac{\partial{\varphi_3}}{\partial{t}}(\sigma,t)\mid{_{s=s_0}}=0$ hold simultaneously.
\end{lemma}

\begin{lemma}
Let $s_0\in\Omega$, if $s_0$ is not a zero of $W_3^{'}(s)$,
then $\frac{\partial{\varphi_3}}{\partial{\sigma}}(\sigma,t)$ and
$\frac{\partial{\varphi_3}}{\partial{t}}(\sigma,t)$ cannot simultaneously equal to 0.
\end{lemma}

The lemma can also be written as: let $s_0\in\Omega$, and it is not a zero of $W_3^{'}(s)$,
if $\frac{\partial{\varphi_3}}{\partial{\sigma}}(\sigma,t)\mid{_{s=s_0}}=0$,
then, $\frac{\partial{\varphi_3}}{\partial{t}}(\sigma,t)\mid{_{s=s_0}}\neq 0$; or,
if $\frac{\partial{\varphi_3}}{\partial{t}}(\sigma,t)\mid{_{s=s_0}}=0$,
then, $\frac{\partial{\varphi_3}}{\partial{\sigma}}(\sigma,t)\mid{_{s=s_0}}\neq 0$.

Because $W_3(s)=iW(s)$, $W_3^{'}(s)=iW^{'}(s)$. Obviously,
two functions $W_3(s)$ and $W(s)$ have the same derivative zeros. In front, we prove the two partial derivatives of argument functions of two functions ($W_3(s)$ and $W(s)$) are separately same.
According to the derivative function relationship which is given in this paragraph,
and according to Lemma 1.7 and Lemma 1.9, we obtain Lemma 1.11.

Let $\Theta=\{s=\sigma+it|s\in\mathbb{C}\cup\{\infty\}$, $
s$ is not zero or pole of $W(s)$, and there is $u(\sigma,t)=0\}$.

\begin{lemma}
Let $s_0\in\Theta$, then $s_0=\sigma_0+it_0$ is the zero of $W^{'}(s)$
if and only if there are $\frac{\partial{\varphi}}{\partial{\sigma}}(\sigma,t)\mid{_{s=s_0}}=0$
and $\frac{\partial{\varphi}}{\partial{t}}(\sigma,t)\mid{_{s=s_0}}=0$.
\end{lemma}

According to the Cauchy-Riemann Equations, apply Lemma 1.7 and Lemma 1.10, we obtain Lemma 1.12.

\begin{lemma}
Let $s_0=\sigma_0+it_0\in\Theta$, if $s_0$ is not zero of $W^{'}(s)$,
then $\frac{\partial{\varphi}}{\partial{\sigma}}(\sigma,t)$ and
$\frac{\partial{\varphi}}{\partial{t}}(\sigma,t)$ cannot simultaneously equal to 0.
\end{lemma}

The lemma can also be written as: Let $s_0=\sigma_0+it_0\in\Theta$, if $s_0$ is not zero of $W^{'}(s)$,
then if $\frac{\partial{\varphi}}{\partial{\sigma}}(\sigma,t)\mid{_{s=s_0}}=0$,
then, $\frac{\partial{\varphi}}{\partial{t}}(\sigma,t)\mid{_{s=s_0}}\neq 0$;
or, if $\frac{\partial{\varphi}}{\partial{t}}(\sigma,t)\mid{_{s=s_0}}=0$,
then, $\frac{\partial{\varphi}}{\partial{\sigma}}(\sigma,t)\mid{_{s=s_0}}\neq 0$.

Let $\Delta=\{s=\sigma+it|s\in\mathbb{C}\cup\{\infty\}$, $s$ is not zero or pole of $W(s)\}$.
We need to prove the results suitable in the whole $\mathbb{C}\cup\{\infty\}$, except the zero and pole of $W(s)$.
Because $\Phi$ is equivalent to $\Delta\cap\{u(\sigma,t)\neq0\}$, and $\Theta$ is equivalent to $\Delta\cap\{u(\sigma,t)=0\}$,
we have $\Phi\cap\Theta=\emptyset$, $\Phi\subset\Delta$, $\Theta\subset\Delta$, and $\Delta=\Phi\cup\Theta$.

Lemma 1.5 and Lemma 1.6 are both except the points which let $u(\sigma,t)$ of $W(s)$ equals to 0. The results in Lemma 1.5 and Lemma 1.6 are true for the points in $\Phi$.
Lemma 1.11 and Lemma 1.12 give the same results under the condition $u(\sigma,t)=0$. The results in Lemma 1.11 and Lemma 1.12 are true for the points in $\Theta$. According to $\Delta=\Phi\cup\Theta$.
Therefore, we combine Lemma 1.5 and Lemma 1.11, obtain Theorem 1.13.
We combine Lemma 1.6 and Lemma 1.12, obtain Theorem 1.14.

\begin{theorem}
Let $s_0=\sigma_0+it_0\in\mathbb{C}\cup\{\infty\}$, and $s_0$ is not zero or pole of $W(s)$,
then $s_0$ is the zero of $W^{'}(s)$
if and only if there are $\frac{\partial{\varphi}}{\partial{\sigma}}(\sigma,t)\mid{_{s=s_0}}=0$
and $\frac{\partial{\varphi}}{\partial{t}}(\sigma,t)\mid{_{s=s_0}}=0$.
\end{theorem}

\begin{theorem}
Let $s_0=\sigma_0+it_0\in\mathbb{C}\cup\{\infty\}$, and $s_0$ is not zero or pole of $W(s)$,
if $s_0$ is not zero of $W^{'}(s)$,
then $\frac{\partial{\varphi}}{\partial{\sigma}}(\sigma,t)$ and
$\frac{\partial{\varphi}}{\partial{t}}(\sigma,t)$ cannot simultaneously equal to 0.
\end{theorem}

The theorem can also be written as: let $s_0=\sigma_0+it_0\in\Delta$, if $s_0$ is not zero of $W^{'}(s)$,
then if $\frac{\partial{\varphi}}{\partial{\sigma}}(\sigma,t)\mid{_{s=s_0}}=0$,
we have $\frac{\partial{\varphi}}{\partial{t}}(\sigma,t)\mid{_{s=s_0}}\neq 0$;
or, if $\frac{\partial{\varphi}}{\partial{t}}(\sigma,t)\mid{_{s=s_0}}=0$, we have
$\frac{\partial{\varphi}}{\partial{\sigma}}(\sigma,t)\mid{_{s=s_0}}\neq0$.

Theorem 1.13 and Theorem 1.14 are the results of the necessary and sufficient condition. And no matter the $s_0$ is real or complex,
and no matter the order of $W(s)$ is finite or infinite, the two theorems are both true.

All the poles of the Gamma function $\Gamma(s)$ are located at the negative integer number,
and there is no zeros of $\Gamma(s)$ in the extended complex plane $\mathbb{C}\cup\{\infty\}$.
This can be obtained directly from its analytic expression\cite{Stein,Carlo}.
According to the Rolle's theorem in the mathematical analysis, we can easily obtain that for the digamma function,
in each pair of adjacent negative integers, there exists a zero of the digamma function $\Psi(s)$,
namely, there exists a zero of $\Gamma^{'}(s)$. But, Rolle's theorem is a sufficient result. It cannot be used to prove there cannot exist the zeros of derivative of a function on other position in the extend complex plane. Theorem 1.13 can solve such problems. Here, we prove the zero of $\Gamma^{'}(s)$ only distributed between each pair of adjacent negative integers, and there cannot exist the zeros of derivative on the other positions in the extend complex plane. Our result let the study of the zero of $\Gamma^{'}(s)$ more precise and comprehensive.

\begin{theorem}
Let $s\in\mathbb{C}\cup\{\infty\}$,
if $s$ is not on the real axis,
then $s$ can not be the zero of $\Gamma^{'}(s)$ or $\Psi^{'}(s)$.
\end{theorem}
\begin{proof}
The gamma function is defined as:
$\Gamma(s)=\frac{1}{se^{rs}\prod^{+\infty}_{k=1}(1+\frac{s}{k})e^{-\frac{s}{k}}}$,
and digamma function $\Psi(s)$ is defined as:
$\Psi(s)=\frac{\Gamma^{'}(s)}{\Gamma(s)}$\cite{Stein,Carlo}.

The argument of $\Gamma(s)$ is: $\arg(\Gamma(s))=-\arg(s)-r\arg(e^{s})-\sum^{+\infty}_{k=1}(\arg(1+\frac{s}{k})-\arg(e^{-\frac{s}{k}}))$.
Let $s=\sigma+it$, the expression can be written as:
$\arg(\Gamma(\sigma+it))=-\arg(\sigma+it)-rt-\sum^{+\infty}_{k=1}(\arg(1+\frac{\sigma+it}{k})-\frac{t}{k})=$
$-\arg(\frac{t}{\sigma})-rt-\sum^{+\infty}_{k=1}(\arg(\frac{t}{\sigma+k})-\frac{t}{k})$.

According to Lemma 1.3, take the two partial derivatives of the argument of the gamma function, we have:
$\frac{\partial{\arg}}{\partial{\sigma}}(\Gamma(\sigma+it))=$
$-\frac{\partial{\arg}}{\partial{\sigma}}(\frac{t}{\sigma})$
$-\sum^{+\infty}_{k=1}\frac{\partial{\arg}}{\partial{\sigma}}(\frac{t}{\sigma+k})$
$=\frac{t}{t^{2}+\sigma^{2}}+\sum^{+\infty}_{k=1}\frac{t}{(\sigma+k)^{2}+t^{2}}$,
$\frac{\partial{\arg}}{\partial{t}}(\Gamma(\sigma+it))=$
$-\frac{\partial{\arg}}{\partial{t}}(\frac{t}{\sigma})-r-$
$\sum^{+\infty}_{k=1}(\frac{\partial{\arg}}{\partial{t}}(\frac{t}{\sigma+k})-\frac{1}{k})$
$=-\frac{\sigma}{t^{2}+\sigma^{2}}-r-\sum^{+\infty}_{k=1}(\frac{\sigma+k}{(\sigma+k)^{2}+t^{2}}-\frac{1}{k})$.

When $t>0$, no matter which value $\sigma$ takes,
there is $t(\frac{1}{t^{2}+\sigma^{2}}+\sum^{+\infty}_{k=1}\frac{1}{(\sigma+k)^{2}+t^{2}})>0$.
So, when $t>0$, $\frac{\partial{\arg}}{\partial{\sigma}}(\Gamma(\sigma+it))>0$.
According to Theorem 1.14, no matter which value $\frac{\partial{\arg}}{\partial{t}}(\Gamma(\sigma+it))$ takes,
in the upper half of the complex plane, there has no zeros of $\Gamma^{'}(s)$ and $\Psi^{'}(s)$.

$\Gamma(s)$ is a conjugate function with real coefficients,
and its poles are all simple\cite{Stein,Carlo}. So, its derivative is also a conjugate function with real coefficients,
and $\Psi(s)$ is also a conjugate function with real coefficients.
So, if we prove that there is no zero of $\Gamma^{'}(s)$ and $\Psi^{'}(s)$ in the upper half of the complex plane,
then, in the lower half of the complex plane, $\Gamma^{'}(s)$ and $\Psi^{'}(s)$ can not have zeros.
\end{proof}

The following theorem has been proved by Brian Conrey\cite{Brian},
here, we use the new theorem to prove it in a new way.

\begin{theorem}
Let $s\in\mathbb{C}\cup\{\infty\}$, if $s$ is not on the critical strip,
then $s$ can not be zero of $\xi^{'}(s)$.
Assuming that the Riemann hypothesis is true,
then all the zeros of $\xi^{'}(s)$ are on the critical line.
\end{theorem}
\begin{proof}
Hadamard proved that Xi-function is expressed as: $\xi(s)=\frac{1}{2}\prod_{\rho}(1-\frac{s}{\rho})$\cite{Norman}.
The argument of xi-function is: $\arg(\xi(s))=\sum_{\rho}(\arg(s-\rho)$$-\arg(-\rho))=$
$\sum_{\rho}(\arg(\frac{t-t_{\rho}}{\sigma-\sigma_{\rho}}))-$ $\arg(\frac{-t_{\rho}}{-\sigma_{\rho}}))$.

According to Lemma 1.3, we do the two partial derivatives of argument of $\xi(s)$
concerning the the real part variable $\sigma$ and imaginary part variable $t$.
We have: $\frac{\partial}{\partial{\sigma}}\arg(\xi(s))=$
$\sum_{\rho}(\frac{\partial}{\partial{\sigma}}\arg(\frac{t-t_{\rho}}{\sigma-\sigma_{\rho}})
-\frac{\partial}{\partial{\sigma}}\arg(\frac{-t_{\rho}}{-\sigma_{\rho}}))$
$=-\sum_{\rho}\frac{t-t_{\rho}}{(\sigma-\sigma_{\rho})^{2}+(t-t_{\rho})^{2}}$,
$\frac{\partial}{\partial{t}}\arg(\xi(s))=
\sum_{\rho}(\frac{\partial}{\partial{t}}\arg(\frac{t-t_{\rho}}{\sigma-\sigma_{\rho}})
-\frac{\partial}{\partial{t}}\arg(\frac{-t_{\rho}}{-\sigma_{\rho}}))$
$=\sum_{\rho}$$\frac{\sigma-\sigma_{\rho}}{(\sigma-\sigma_{\rho})^{2}+(t-t_{\rho})^{2}}$.

Because all of the non-trivial zeros are all distributed outside the critical strip,
namely, for all of the points which are distributed outside the right side of the critical strip,
$\sigma>\sigma_{\rho}$, so, $\frac{\sigma-\sigma_{\rho}}{(\sigma-\sigma_{\rho})^{2}+(t-t_{\rho})^{2}}>0$,
$\frac{\partial}{\partial{t}}\arg(\xi(s))=
\sum_{\rho}\frac{\sigma-\sigma_{\rho}}{(\sigma-\sigma_{\rho})^{2}+(t-t_{\rho})^{2}}>0$.
According to Theorem 1.14, all the points which distributed on this zone cannot be the zeros of $\xi^{'}(s)$.

For all the points which are distributed outside the left side of the critical strip, there are $\sigma<0$,
$\sigma-\sigma_{\rho}<0$, so, $\frac{\sigma-\sigma_{\rho}}{(\sigma-\sigma_{\rho})^{2}+(t-t_{\rho})^{2}}<0$,
$\frac{\partial}{\partial{t}}\arg(\xi(s))=
\sum_{\rho}\frac{\sigma-\sigma_{\rho}}{(\sigma-\sigma_{\rho})^{2}+(t-t_{\rho})^{2}}<0$.
According to Theorem 1.14, all of points which distributed on this zone cannot be the zeros of $\xi^{'}(s)$.

Sum up the above results, we can obtain the result, there are no zeros of $\xi^{'}(s)$ outside the critical strip.

If all of the non-trivial zeros are on the critical line,
for all points which are distributed at the right side of the critical line, there are $\sigma>\frac{1}{2}$,
so, $\frac{\sigma-\frac{1}{2}}{(\sigma-\frac{1}{2})^{2}+(t-t_{\rho})^{2}}>0$,
$\frac{\partial}{\partial{t}}\arg(\xi(s))=
\sum_{\rho}\frac{\sigma-\frac{1}{2}}{(\sigma-\frac{1}{2})^{2}+(t-t_{\rho})^{2}}>0$.
According to Theorem 1.14, all of points which distributed on this zone cannot be the zeros of $\xi^{'}(s)$.

For all the points which are distributed at the left side of the critical line,
there are $\sigma<\frac{1}{2}$, so, $\frac{\sigma-\frac{1}{2}}{(\sigma-\frac{1}{2})^{2}+(t-t_{\rho})^{2}}<0$,
$\frac{\partial}{\partial{t}}\arg(\xi(s))=
\sum_{\rho}\frac{\sigma-\frac{1}{2}}{(\sigma-\frac{1}{2})^{2}+(t-t_{\rho})^{2}}<0$.
According to Theorem 1.14, all of points which distributed on this zone cannot be the zeros of $\xi^{'}(s)$.
So, assuming that the Riemann hypothesis is true, then all the zeros of $\xi^{'}(s)$ are on the critical line.
\end{proof}

\bibliographystyle{unsrt}

\begin{thebibliography}{99}

\bibitem{Fischer} W. Fischer and I. Lieb. A course in complex analysis. From Basic Results
to Advanced Topics. Translated by Cannizzo, J. Vieweg+Teubner Verlag, (2012)

\bibitem {Hayman1} W. K. Hayman, Picard Values of Meromorphic Functions and Their Derivatives. Annals of Mathematics, 70, 9-42 (1959)

\bibitem {Edrei2} A. Edrei, Meromorphic functions with three radially distributed values, Trans. Amer. Math. Soc. 78, 276-293 (1955)

\bibitem {Ostrovskii} A. A. Gol’dberg and I. V. Ostrovskii, Value Distribution of Meromorphic Functions, Nauka,
Moscow, 1970 (Russian); English transl., Amer. Math. Soc. Providence, RI, (2008)

\bibitem {Edrei1}  A. Edrei and W. H. J. Fuchs, On meromorphic functions with regions free of poles and zeros, Acta
Math. 108 113-145 (1962)

\bibitem {Hayman2} W. K. Hayman, Meromorphic Functions, Clarendon Press, Oxford, (1964)

\bibitem {Ki} H. Ki and Y.-O. Kim, On the number of nonreal zeros of real entire functions and the Fourier-
Pólya conjecture, Duke Math. J. 104 45-73 (2000)

\bibitem {Langley1} J. K. Langley, Derivatives of meromorphic functions of finite order, Comput. Methods Funct.
Theory 14 195-207 (2014)

\bibitem {Langley2} J. K. Langley, Non-real zeros of derivatives of meromorphic functions, J. Anal. Math.
133 183-228 (2017)

\bibitem {Konstantin} Konstantin M. Dyakonov, Zeros of analytic functions, with or without multiplicities,
Math. Ann. 352:625-641 (2012)

\bibitem {Yamanoi} Katsutoshi Yamanoi, Zeros of higher derivatives of meromorphic functions in the
complex plane. Proc. London Math. Soc. (3) 106 703-780 (2013)

\bibitem {Xuecheng} Xuecheng Pang, Shahar Nevo and Lawrence Zalcman, Derivatives of Meromorphic Functions
with Multiple Zeros and Rational Functions. Comput. Methods Funct. Theory, 8(2), 483-491 (2008)

\bibitem {Lars} Lars Ahfors, Complex Analysis: An Introduction To The Theory Of Annalytic Function Of One Complex Variable, MeGraw-Hill (1979)

\bibitem {Stein} E.M. Stein, R. Shakarchi, Complex analysis, Princeton Lectures in Analysis, II,
Princeton University Press, Princeton, NJ, (2003)

\bibitem {Carlo} Carlo Viola, An Introduction to Special Functions, Springer International Publishing, (2016)

\bibitem {Norman} Norman Levinson, Zeros of derivative of Riemann's $\xi(s)$ function,
Bull. Amer. Math. Soc. 80 (5)951-954 (1974)

\bibitem {Brian} Brian Conrey, Zeros of derivatives of Riemann’s xi-function on the critical line,
J. Number Theory 16 no.1, 49-74 (1983)


\end{thebibliography}

\end{document}